\documentclass[11pt]{article}
 
\usepackage{times}
\usepackage{amssymb}
\usepackage{amsmath}
\usepackage{amsbsy}
\usepackage{amsfonts}
\usepackage{eucal}
\newtheorem{theo}{Theorem}[section]
\newtheorem{lm}[theo]{Lemma}
 
\newcommand{\cc}{{\mathbb C}}
\newcommand{\ket}[1]{\left | #1 \right \rangle}
\newcommand{\wt}[1]{{\mathsf wt}(#1)}
\newcommand{\bsym}[1]{{\boldsymbol #1}}

\newcommand{\qed}{\mbox{\rule{1.6mm}{4.3mm}}}

\begin{document}
\title{\bf Bounds for weight distribution of weakly self--dual codes\thanks{This work 
was supported in part by grants from the Revolutionary Computing
group at JPL (contract \#961360), and from the DARPA Ultra program
(subcontract from Purdue University \#530--1415--01).}}

\author{ 
Vwani P. Roychowdhury\thanks{V. Roychowdhury is with the Electrical Engineering Department, 
UCLA, Los Angeles, CA 90095 (e-mail: vwani@ee.ucla.edu).} and
Farrokh Vatan\thanks{F. Vatan is with the Electrical Engineering Department, UCLA, 
Los Angeles, CA 90095 (e-mail: vatan@ee.ucla.edu).}}

%\markboth{IEEE Transactions On Information Theory, Vol. XX, No. Y, Month 2000}
%{Vatan and Roychowdhury: Bounds for weight distribution of weakly self--dual codes}

\maketitle

\begin{abstract}
Upper bounds are given for the weight distribution of binary weakly self--dual codes.
To get these new bounds, we introduce a novel method of utilizing unitary operations
on Hilbert spaces. This method is motivated by recent progress on quantum computing.
This new approach leads to much simpler proofs for such genre of bounds on the 
weight distributions of certain classes of codes. Moreover, in some cases, our bounds are 
improvements on the earlier bounds. These improvements are achieved, either by 
extending the range of the weights over which the bounds apply, or by extending the class of 
codes subjected to these bounds. 

\end{abstract}

\section{Background}

For a random linear $[n,k]$ code $\cal C$ with weight distribution 
$(A_0,A_1,\ldots,A_n)$ it is known that the expected value of the normalized weight
distribution, i.e., $\frac{1}{2^k}A_w$, is the same as the normalized binomial
distribution $\frac{1}{2^n}\binom{n}{w}$ (see, e.g., [14, p. 287]).
So for such a code the expected value for the number of codewords of weight $w$ is
$\frac{1}{2^{(n-k)}}\binom{n}{w}$. The problem to determine which explicit
classes of codes have binomial weight distribution has been investigated in several papers
(see references). Most of these results are about the BCH codes, their
extensions, or their dual codes. For example, it is shown that for these codes the number
of codewords of weight $w$ is 
\[ A_w=\frac{1}{2^{n-k}}\left(\binom{n}{w}+E_w \right), \]
for $w$ in some range, and the error term $E_w$ tends to zero when $n$ tends to infinity.  

There are also bounds for other classes of liner codes.
Let $\cal C$ be a doubly--even self--dual $[n,n/2,d]$ 
code $\cal C$ with weight distribution $(A_0,A_1,\ldots,A_n)$. 
Let $\delta=\frac{d}{n}$, then in [11] it is shown that
\begin{equation}
 A_w \leq 2^{\left(H_2(\frac{w}{n})-\frac{1}{2}\right)n}, 
\label{bound1}
\end{equation}
if $\frac{w}{n}\in[c,1-c]$, where  
\begin{equation}
 c=\frac{1}{2}-\sqrt{\frac{6\delta-1+\sqrt{1-8\delta+32\delta^2}}{8(1-\delta)}}. 
\label{cond-c}
\end{equation}
This shows that for this class of codes, the weight distribution around $\frac{n}{2}$
is upper--bounded by the binomial distribution.

In [12] an upper bound for weight distribution of the {\em dual} of
extended BCH codes is given. Let $\cal C$ of length $n=2^m$ be the dual of the extended 
code of a $t$--error correcting BCH code. Then $\dim({\cal C})\leq mt$; i.e., 
$\vert{\cal C}\vert\leq n^t$. If $(A_0,A_1,\ldots,A_n)$ is the weight distribution 
of $\cal C$, then for $w>\sqrt{\frac{n}{t+1}}+2$ we have
\begin{equation}
\textstyle A_w\leq \frac{4\sqrt{2\pi n}n^t}{\vert 2\sqrt{n(t+1)}-n+2w\vert} 
    \cdot e^{-\frac{(n-2w)^2}{n}}\left ( 1+O\left(\frac{1}{t}\right)\right ) . 
\label{bound2}
\end{equation}

In [18] several bounds for the weight distribution of subfield subcodes of
algebraic--geometric codes is given. This class of codes contains important
classes of codes, such as the binary BCH and Goppa codes. For an $[n,k]$ binary code
of this type, with weight distribution $(A_0,A_1,\ldots,A_n)$, they derived the following
bound
\begin{equation}
 \left\vert A_w-\frac{1}{2^{n-k}}\binom{n}{w}\right\vert < c_1 n^{\frac{w}{2}}, 
\label{bound3}
\end{equation}
and for the special case of $w=\frac{n}{2}$, they get the following bound
\begin{equation}
 \left\vert A_{\frac{n}{2}}-\frac{1}{2^{n-k}}\binom{n}{\frac{n}{2}}\right\vert 
             < c_2{\binom{n}{\frac{n}{2}}}^{\frac{1}{2}}n^{\frac{1}{4}}. 
\label{bound4}
\end{equation}
Here $c_1$ and $c_2$ are two constants, both much larger than $\sqrt{e}$.

In this paper, {\it we apply a novel approach based on unitary operations on 
Hilbert spaces}, 
and {\it derive bounds for the weight distribution of another class of
linear codes}. In particular, we study the class of weakly self--dual codes, 
i.e., the class of
codes $\cal C$ such that ${\cal C}\subseteq{\cal C}^\perp$. We show that, for 
$0<w<\frac{n}{2}$, the number $A_w$ of codewords of weight $w$ in $\cal C$ satisfies
the following bounds
\begin{equation}
 A_w \leq 2^{\frac{1}{2}H_2(\frac{w}{n})n},
\label{bound5}
\end{equation}
and
\begin{equation}
  A_w \leq \sqrt{e}(n-w+1)^{\frac{w}{2}}.
\label{bound6}
\end{equation}

If we compare our bounds with previously known bounds (1), (3),
and (4), we realize that these new bounds, for some values of $w$
and in the intersection of their corresponding classes of codes, 
give a better estimate than the old bounds. For example, (6) 
holds for any value of $w$, $0<w<\frac{n}{2}$, and it applies also to 
the special case of doubly--even self--dual codes, and, hence, comparable to the 
bounds in [11]. The bound in (1) applies
to doubly--even self--dual codes, but holds only in the interval
$[c,1-c]$, with $c$ defined by (2). So if 
$\frac{d}{n} \leq {H_2}^{-1}(\frac{1}{2})=0.1100\cdots$ then $c>0.27$. Hence, for this 
choice of $\frac{d}{n}$, 
(6) is a better bound than (1) for $\frac{w}{n}<0.27$, since
(1) does not even hold for these values of $w$; we should note, however, that 
(1) is a better bound if $\frac{w}{n}>0.27$.
One can also verify that (7) gives a better bound than (4),
as the constant $c_1$ is much larger than $\sqrt{e}$.

\section{Unitary operations on Hilbert spaces}

We derive our bounds via actions of unitary operations on Hilbert spaces.
Toward this end, we find it more comfortable to use the language of ``bra--ket''
of quantum mechanics (see, e.g., [3]) as it is used in the theory of
quantum computation (see, e.g., [15] and [1]). 
We briefly describe the necessary notions and notations.

Consider the two--dimensional Hilbert space ${\cal H}=\cc^2$. 
We denote its standard basis by $\{\ket{0},\ket{1}\}$; i.e., $\ket{0}=(1,0)$ 
and $\ket{1}=(0,1)$. We consider also the tensor product 
\[ {\cal H}^{\otimes n}=\cc^2\otimes\cc^2\otimes\cdots\otimes\cc^2 \]
of $n$ copies of $\cc^2$. So ${\cal H}^{\otimes n}$ is a $2^n$--dimensional Hilbert 
space isomorphic with $\cc^{2^n}$.
We represent the standard basis of ${\cal H}^{\otimes n}$ by $2^n$ products of
the form 
\[ \ket{c_1}\otimes\ket{c_2}\otimes\cdots\otimes\ket{c_n}, \]
where $c_i\in\{0,1\}$. For simplicity, we write $\ket{c_1c_2\cdots c_n}$ instead of
$\ket{c_1}\otimes\ket{c_2}\otimes\cdots\otimes\ket{c_n}$. The Euclidean length 
on ${\cal H}^{\otimes n}$ is defined in the natural way; we denote the length of
the vector $\ket{a}\in{\cal H}^{\otimes n}$ by $\Vert\ket{a}\Vert$. For example, 
${\cal H}^{\otimes 2}=\cc^2\otimes\cc^2$ is a 4--dimensional Hilbert space with
\[ \{\ket{00},\ket{01},\ket{10},\ket{11}\} \]
as its standard basis. If 
\[ \ket{a}=a\ket{00}+b\ket{01}+c\ket{10}+d\ket{11}\in{\cal H}^{\otimes 2} \]
then $\Vert\ket{a}\Vert=\left(aa^*+bb^*+cc^*+dd^*\right)^{\frac{1}{2}}$,
where $^*$ stands for the complex conjugate.

An $m\times m$ matrix $M$ is {\em unitary} if $M^\dagger\cdot M=I_m$,
where $I_m$ is the identity matrix and $M^\dagger$ is the {\em adjoint} 
matrix of $M$; i.e., $M^\dagger=\left(M^{\mathrm tr}\right)^*$, where
``tr'' denotes the transpose.  
A linear operation on the $m$--dimensional Hilbert space 
$\cc^m$ is called a {\em unitary operation} if it is represented by a unitary
matrix. Note that every unitary operation is a {\em length--preserving} operation.

We denote the group of $m\times m$ unitary matrices by
$\mbox{\bf U}(m)$. The general form of a matrix in $\mbox{\bf U}(2)$ is
\[\begin{pmatrix}
  e^{i(\alpha+\gamma)}\cos\theta & e^{i(\beta+\gamma)}\sin\theta \\
 -e^{-i(\beta-\gamma)}\sin\theta & e^{-i(\alpha-\gamma)}\cos\theta
  \end{pmatrix} . \]
Suppose that $U_1,U_2\in\mbox{\bf U}(2)$, then the tensor product
$U_1\otimes U_2$ in $\cc^2\otimes\cc^2$ is defined in the natural way. For example, 
if $U_1(\ket{0})=a\ket{0}+b\ket{1}$ and $U_2(\ket{1})=c\ket{0}+d\ket{1}$ then
\begin{eqnarray*}
 U_1\otimes U_2(\ket{01}) & = & U_1(\ket{0})\otimes U_2(\ket{1}) \\
  & = & (a\ket{0}+b\ket{1})\otimes(c\ket{0}+d\ket{1}) \\
  & = & ac\ket{00}+ad\ket{01}+bc\ket{01}+bd\ket{11}.
\end{eqnarray*}

\section{Bounds for weights}

For any real number $\theta$, consider the unitary operation 
$R_\theta\in\mbox{\bf U}(2)$ defined by the following matrix
\[ R_\theta=\begin{pmatrix} \sin\theta & \cos\theta \\ \cos\theta & -\sin\theta
            \end{pmatrix} . \]
Then the action of $R_\theta$ on a vector $\ket{c}$, $c\in\{0,1\}$, can be
written as follows:
\[ R_\theta(\ket{c})=\sum_{a\in\{0,1\}}(-1)^{ac}({\sin\theta})^{1-c-a+2ac}
                      	({\cos\theta})^{c+a-2ac}\ket{a} . \]
Now let $\bsym{S}_\theta={R_\theta}^{\otimes n}$, i.e., the tensor product of
$n$ copies of $R_\theta$. Thus $\bsym{S}_\theta\in\mbox{\bf U}(2^n)$.
The following lemma provides a closed form for the action of $\bsym{S}_\theta$.

\begin{lm}
For any vector $\ket{\bsym{c}}$ in the standard basis, i.e., $\bsym{c}\in\{0,1\}^n$, 
we have
\begin{equation}
 \bsym{S}_\theta(\ket{\bsym{c}})= 
 \sum_{\bsym{a}\in\{0,1\}^n}(-1)^{\bsym{c}\cdot\bsym{a}}
 (\sin\theta)^{n-\wt{\bsym{c}+\bsym{a}}} 
 (\cos\theta)^{\wt{\bsym{c}+\bsym{a}}} \ket{\bsym{a}},
\label{rt}
\end{equation}
where $\bsym{c}\cdot\bsym{a}$ is the inner product of $\bsym{c}$ and $\bsym{a}$ as real
vectors, $\wt{\bsym{x}}$ denotes the Hamming weight of a binary vector $\bsym{x}$,
and the addition of binary vectors is considered over GF(2). 
\end{lm}

{\bf Proof.}
We first show that
\begin{equation}
 \bsym{S}_\theta(\ket{\bsym{c}})= 
  \sum_{\bsym{a}\in\{0,1\}^n}(-1)^{\bsym{c}\cdot\bsym{a}}
  (\sin\theta)^{n-\wt{\bsym{c}}-\wt{\bsym{a}}+2\bsym{c}\cdot\bsym{a}}
  (\cos\theta)^{\wt{\bsym{c}}+\wt{\bsym{a}}-2\bsym{c}\cdot\bsym{a}} \ket{\bsym{a}} ,
\label{rt1}
\end{equation}
We prove the identity (9) only for $n=2$; the proof in the general case is
quite similar. Let $\bsym{c}=(c_1,c_2)$ and $\bsym{a}=(a_1,a_2)$.
\begin{eqnarray*}
 \bsym{S}_\theta(\ket{\bsym{c}}) & = & 
        R_\theta(\ket{c_1})\otimes R_\theta(\ket{c_2})  \\
 & = & \Big(\sum_{a_1\in\{0,1\}}(-1)^{a_1c_1}({\sin\theta})^{1-c_1-a_1+2a_1c_1} 
      ({\cos\theta})^{c_1+a_1-2a_1c_1}\ket{a_1}\Big)\otimes \\
 & &   \Big(\sum_{a_2\in\{0,1\}}(-1)^{a_2c_2}({\sin\theta})^{1-c_2-a_2+2a_2c_2} 
      ({\cos\theta})^{c_2+a_2-2a_2c_2}\ket{a_2}\Big) \\
 & = & \sum_{a_1,a_2\in\{0,1\}}(-1)^{a_1c_1+a_2c_2} 
      ({\sin\theta})^{2-(c_1+c_2)-(a_1+a_2)+2(a_1c_1+a_2c_2)} \\
 & &   ({\cos\theta})^{(c_1+c_2)+(a_1+a_2)-2(a_1c_1+a_2c_2)}\ket{a_1a_2}, 
\end{eqnarray*} 
which is the of form of (9). We note that
\[ \wt{\bsym{c}}+\wt{\bsym{a}}-2\bsym{c}\cdot\bsym{a}=\wt{\bsym{c}+\bsym{a}}. \]
Therefore, (9) can be rewritten as (8). \qed

Now we consider a linear code ${\cal C}\subseteq\{0,1\}^n$ of dimension $k$,
and the corresponding {\em unit} vector
\[ \ket{\cal C}=\frac{1}{\sqrt{2^k}}\sum_{\bsym{c}\in{\cal C}}\ket{\bsym{c}} .\]
By applying the unitary operation $\bsym{S}_\theta$ on the unit vector $\ket{\cal C}$
we get
\[ \bsym{S}_\theta(\ket{\cal C})= 
     \frac{1}{\sqrt{2^k}}\sum_{\bsym{c}\in{\cal C}}
     \sum_{\bsym{a}\in\{0,1\}^n}(-1)^{\bsym{c}\cdot\bsym{a}} 
     (\sin\theta)^{n-\wt{\bsym{c}+\bsym{a}}} 
     (\cos\theta)^{\wt{\bsym{c}+\bsym{a}}} \ket{\bsym{a}}. \]
This can be rewritten as follows
\[   \bsym{S}_\theta(\ket{\cal C})= 
     \frac{1}{\sqrt{2^k}}\sum_{\bsym{a}\in{\cal C}^\perp}
     \left ( \sum_{\bsym{c}\in{\cal C}}
     (\sin\theta)^{n-\wt{\bsym{c}+\bsym{a}}} 
     (\cos\theta)^{\wt{\bsym{c}+\bsym{a}}} \right ) \ket{\bsym{a}} 
     +\ket{\mathrm{remainder}}.\]
Since $\left\Vert\bsym{S}_\theta(\ket{\cal C})\right\Vert =1$, the following
lemma follows.

\begin{lm}
For any $k$--dimensional linear code $\cal C$ of length $n$, and any real number $\theta$
we have
\begin{equation}
  \frac{1}{2^k}\sum_{\bsym{a}\in{\cal C}^\perp}
     \left ( \sum_{\bsym{c}\in{\cal C}}
     (\sin\theta)^{n-\wt{\bsym{c}+\bsym{a}}} 
     (\cos\theta)^{\wt{\bsym{c}+\bsym{a}}} \right )^2\leq 1.
\label{ineq1}
\end{equation}
\label{lm1}
\end{lm}

\begin{lm}
For any $k$--dimensional weakly self--dual linear code $\cal C$ of length $n$, i.e.,  
${\cal C}\subseteq {\cal C}^\perp$, and any real number $\theta$ we have
\begin{equation}
     \left\vert\sum_{\bsym{c}\in{\cal C}^\perp}
     (\sin\theta)^{n-\wt{\bsym{c}}} 
     (\cos\theta)^{\wt{\bsym{c}}}\right\vert\leq 2^{(n-2k)/2}.
\label{ineq2}
\end{equation}
\label{lm2}
\end{lm}

{\bf Proof.}
We apply Lemma III.2 to the code ${\cal C}^\perp$. The result is
\[ \frac{1}{2^{n-k}}\sum_{\bsym{a}\in{\cal C}}
   \left ( \sum_{\bsym{c}\in{\cal C}^\perp}(\sin\theta)^{n-\wt{\bsym{c}+\bsym{a}}} 
   (\cos\theta)^{\wt{\bsym{c}+\bsym{a}}} \right )^2\leq 1. \]
Since ${\cal C}\subseteq{\cal C}^\perp$, for every $\bsym{a}\in{\cal C}$ we have
\[  \sum_{\bsym{c}\in{\cal C}^\perp}(\sin\theta)^{n-\wt{\bsym{c}+\bsym{a}}} 
     (\cos\theta)^{\wt{\bsym{c}+\bsym{a}}}= 
     \sum_{\bsym{c}\in{\cal C}^\perp}
     (\sin\theta)^{n-\wt{\bsym{c}}} 
     (\cos\theta)^{\wt{\bsym{c}}}, \]
and the lemma follows. \qed

\begin{lm}
Let $\cal C$ be a weakly self--dual code with weight distribution 
$(A_0,A_1,\ldots,A_n)$. Then for any $0<\lambda<1$ we have
\begin{equation}
 \sum_{j=0}^{n/2}A_{2j}\lambda^j\leq (1+\lambda)^{n/2}, \qquad 0<\lambda<1.
\label{ineq6}
\end{equation}
\end{lm}

{\bf Proof.}
We first apply the following form of the MacWilliams identity (see [14])
to get a handier form of inequality (11):
\[ \sum_{\bsym{u}\in{\cal C}^\perp}x^{n-\wt{\bsym{u}}}y^{\wt{\bsym{u}}}=
   \frac{1}{\vert{\cal C}\vert}
   \sum_{\bsym{u}\in{\cal C}}(x+y)^{n-\wt{\bsym{u}}}(x-y)^{\wt{\bsym{u}}} .\]
This way, from (11), we obtain the following inequality:
\begin{equation}
     \left\vert\sum_{\bsym{c}\in{\cal C}}
     (\sin\theta+\cos\theta)^{n-\wt{\bsym{c}}} 
     (\sin\theta-\cos\theta)^{\wt{\bsym{c}}}\right\vert\leq 2^{n/2}.
\label{ineq3}
\end{equation}

Now suppose that $\frac{\pi}{4}<\theta<\frac{\pi}{2}$, 
$\sqrt{u}=\sin\theta+\cos\theta$, and 
$\sqrt{v}=\sin\theta-\cos\theta$. Thus $1<u<2$ and $v=2-u$. 
Then inequality (13) can be written as
\begin{equation}
     \sum_{\bsym{c}\in{\cal C}}
     {\sqrt{\smash[b]{u}}\,}^{n-\wt{\bsym{c}}} 
     {\sqrt{\smash[b]{2-u}}\,}^{\wt{\bsym{c}}}\leq 2^{n/2}.
\label{ineq4}
\end{equation}
Since $\cal C$ is weakly self--dual then $A_j=0$, for every odd index $j$. 
Therefore, (14) can be written as follows (assume $n$ is even):
\begin{equation}
 \sum_{j=0}^{n/2}A_{2j}u^{\frac{n}{2}-j}(2-u)^j\leq 2^{n/2},\qquad 1<u<2.
\label{ineq5}
\end{equation}
Let $\frac{2-u}{u}=\lambda$, then $0<\lambda<1$ and (15) became
\[ \sum_{j=0}^{n/2}A_{2j}\lambda^j\leq (1+\lambda)^{n/2}. \ \qed \]

\begin{theo}
For every weakly self--dual code $\cal C$ with weight distribution 
$(A_0,A_1,\ldots,A_n)$ we have
\begin{equation}
 A_w \leq 2^{\frac{1}{2}H_2(\frac{w}{n})n},\qquad 0<w<\frac{n}{2} .
\label{b1}
\end{equation}
and
\begin{equation}
 A_w \leq \sqrt{e}(n-w+1)^{w/2},\qquad 0<w<\frac{n}{2} .
\label{b2}
\end{equation}
\end{theo}

{\bf Proof.}
For $0<w<\frac{n}{2}$, let $\alpha=\frac{w}{n}$ and $A_w=2^{\beta n}$. 
From (12) it follows
\[ 2^{\beta n}\lambda^{\frac{w}{2}}\leq (1+\lambda)^{\frac{n}{2}} .\]
Therefore,
\begin{equation}
  \beta\leq -\frac{1}{2}\alpha\log_2\lambda+\frac{1}{2}\log_2(1+\lambda),\qquad
  0<\lambda<1.
\label{ineq7}
\end{equation}

For fixed $0<\alpha<\frac{1}{2}$, let 
\[ F(\lambda)=-\frac{1}{2}\alpha\log_2\lambda+\frac{1}{2}\log_2(1+\lambda). \]
A simple calculation shows that, on the interval $0<\lambda<1$,
the minimum of $F(\lambda)$ is achieved for $\lambda=\frac{\alpha}{1-\alpha}$.
Thus, for $0<\alpha<\frac{1}{2}$,
\[ \beta\leq\frac{1}{2}\left ( -\alpha\log_2\alpha-(1-\alpha)
   \log_2(1-\alpha)\right ) = \frac{1}{2}H_2(\alpha). \]
Therefore, we have proved (16).

On the other hand, from (12) we have
\[ A_w\leq \left ( 1+\frac{1}{\lambda}\right )^{w/2}(1+\lambda)^{(n-w)/2} .\]
Let $\lambda=\frac{1}{n-w}$, and note that $\left ( 1+\frac{1}{t} \right )^t\leq e$.
This proves (17). \qed

%\bibliography{weightbib}

\end{document}